\documentstyle[12pt]{article}
\input amssym.def
\input amssym.tex
\parskip=\medskipamount
\newtheorem{definition}{Definition}[section]
\newtheorem{lemma}[definition]{Lemma}
\newtheorem{theorem}[definition]{Theorem}
\newtheorem{proposition}[definition]{Proposition}
\newtheorem{corollary}[definition]{Corollary}
\newtheorem{remark}[definition]{Remark}

\newtheorem{comment}[definition]{Comment}
\newtheorem{claim}[definition]{Claim}
\newtheorem{convention}[definition]{Convention}

 1   %Caracteres g"ticos.
\font\ddpp=msbm10  scaled \magstep 1  %Caracteres "doble palo".

\newenvironment{proof}{\noindent{\bf Proof~:}}{\QED\medskip}
\def\QED{\hskip0.1em\hfill\null\ \null\nobreak\hfill
\kern3pt\lower1.8pt\vbox{\hrule\hbox   {\vrule\kern1pt\vbox{\kern1.7pt
\hbox{$\scriptstyle QED$}\kern0.2pt}\kern1pt\vrule}\hrule}}
\def\QQED{\fbox
{$\scriptstyle QED$\ }}
\def\R{\hbox{\ddpp R}}               %Numeros reales
\def\T{\hbox{\ddpp T}}               %toros
\def\L{\hbox{\ddpp L}}
\def\eps{\varepsilon}

\def\Im{{\rm Im}}
\def\Ker{{\rm Ker}}
\def\dim{{\rm dim}}
\def\hr{{\rm hr}}
\def\rk{{\rm rank}}

\def\hqed{\hfill\hfill$\square$}
\def\osi{\Omega_{{\rm sympl}}}

\def\diagram#1{\def\normalbaselines{\baselineskip=0pt
\lineskip=10pt\lineskiplimit=1pt}    \matrix{#1}}

\def\hfl#1#2{\smash{\mathop{\hbox to 12 mm{\rightarrowfill}}
\limits^{\scriptstyle#1}_{\scriptstyle#2}}}

\def\vvfl#1#2{\llap{$\scriptstyle #1$}\left\downarrow
\vbox to 6 mm{}\right.\rlap{$\scriptstyle #2$}}

\newcommand{\ba}{\begin{array}}
\newcommand{\ea}{\end{array}}
\newcommand{\hb}[1]{\hbox{\normalshape #1}}
\def\qq{\qquad}

\newcommand{\op}{\oplus}

\newcommand{\SS}{{\cal S}}

\newcommand{\g}{{\frak g}}

\begin{document}
\baselineskip=.55cm
\title{On Symplectically Harmonic Forms on Six-dimensional Nilmanifolds}
\author{R. IB\'A\~NEZ $^{1}$, YU. RUDYAK $^{2}$, A. TRALLE $^{3}$,
L. UGARTE $^{4}$
\\[10pt]
{\small \it $^1$Departamento de Matem\'aticas, Facultad de
Ciencias,}\\[-8pt]
{\small \it Universidad del Pais Vasco,}\\[-8pt]
{\small \it Apdo. 644, 48080 Bilbao, Spain,}\\[-8pt]
{\small \it E-mail:  mtpibtor@lg.ehu.es}\\[8pt]
{\small \it $^2$ FB6/ Mathematik,} \\[-8pt]
{\small \it Universit\"at Siegen,}\\[-8pt]
{\small \it 57068 Siegen, Germany,}\\[-8pt]
{\small \it E-mail: rudyak@mathematik.uni-siegen.de}\\[8pt]
{\small\it $^3$Department of Mathematics,
}\\[-8pt]
{\small\it University of Olsztyn,}\\[-8pt]
{\small\it 10561 Olsztyn, Poland,}\\[-8pt]
{\small\it E-mail: tralle@tufi.wsp.olsztyn.pl}\\[8pt]
{\small \it $^4$Departamento de Matem\'aticas, Facultad de
Ciencias,}\\[-8pt]
{\small\it Universidad de Zaragoza,}\\[-8pt]
{\small\it 50009 Zaragoza, Spain,}\\[-8pt]
{\small\it E-mail: ugarte@posta.unizar.es}
}
\date{\empty}

\maketitle

\begin{abstract}

In the present paper we study the variation of the dimensions $h_k$ of
spaces of symplectically harmonic cohomology classes
(in the sense of Brylinski) on closed symplectic manifolds. We give a
description of such variation for all 6-dimensional nilmanifolds
equipped with symplectic forms. In particular, it turns out that
certain 6-dimensional nilmanifolds possess families of homogeneous
symplectic forms $\omega_t$ for which numbers $h_k(M,\omega_t)$ vary
with respect to $t$. This gives an affirmative answer to a question
raised by Boris Khesin and Dusa McDuff. Our result is in contrast with
the case of 4-dimensional nilmanifolds which do not admit such
variations by a remark of Dong Yan. 
\end{abstract}

\begin{quote}
{\it Mathematics Subject Classification} (1991): 53C15.

{\it PACS numbers}: 02.40.Ma, 03.20.+i, 0.3.65.-w

{\it Key words and phrases}: symplectically harmonic form,
Brylinski's conjecture, Lefschetz map
\end{quote}

\newpage

\setcounter{section}{0}
\section{Introduction}
\setcounter{equation}{0}

Given a symplectic manifold $(M^{2m},\omega)$, we denote by $[\omega]\in
H^2(M)$ the de Rham cohomology class of $\omega$. Here the notation
$M^{2m}$ means that $M$ is a $2m$-dimensional manifold.
Furthermore, we denote by $L_{\omega}:\Omega^k(M) \to
\Omega^{k+2}(M)$ the multiplication by $\omega$ and by
$L_{[\omega]}:H^{k}(M) \to H^{k+2}(M)$ the induced homomorphism in the
de Rham cohomology $H^*(M)$ of $M$. As usual, we write $L$ instead of
$L_{\omega}$ or $L_{[\omega]}$ if there is no danger of confusion. We
say that a symplectic manifold $(M^{2m},\omega)$ satisfies the {\it
Hard Lefschetz condition} if, for every $k$, the homomorphism
$$
L^k: H^{m-k}(M) \to H^{m+k}(M)
$$
is surjective. In view of the Poincar\'e duality, for closed
manifolds $M$ it means that every $L^k$ is an isomorphism.

In 1988 J. L. Brylinski \cite{Br} introduced the concept of
symplectically harmonic forms (resp. Poisson harmonic forms), defined
for any symplectic manifold (resp. Poisson manifold). Further he
conjectured that on compact symplectic manifolds, every de Rham
cohomology class has a symplectically harmonic representative. In
fact, this conjecture asks about the possibility of constructing of a
symplectic Hodge theory.

Brylinski proved that this conjecture is true for compact K\"ahler
manifolds.
 However, it is not true in general, as it was shown in
\cite{FIL,Ma,Yan}. Mathieu \cite{Ma} proved the following theorem,
which applies to an arbitrary (not necessarily compact) symplectic
manifold.

\begin{theorem}
A symplectic manifold $(M^{2m},\omega)$ satisfies the Brylinski
conjecture if and only if it satisfies the Hard Lefschetz condition.
In other words, the following assertions are equivalent:
\par{\rm (i)} every de Rham cohomology class has a symplectically harmonic
representative;
\par {\rm (ii)} for every $k\leq m$, the homomorphism $L^k:
H^{m-k}(M)\longrightarrow H^{m+k}(M)$ is surjective.
\end{theorem}

Mathieu's proof involves the representation theory of quivers and Lie
superalgebras. An alternative and nice proof can be found in the paper
of Yan \cite{Yan}, who studies a special type of infinite dimensional
$\frak s\frak l(2)$-representation and, basing on this, proves a duality theorem for symplectically harmonic forms which, in turn, implies Theorem 1.1.

Given a symplectic manifold $(M^{2m},\omega)$, let
$\Omega^k_{\hr}(M)=\Omega^k_{\hr}(M,\omega)$ denote the subspace of
symplectically harmonic forms of $\Omega^k(M)$. The form $\alpha$ is
called symplectically harmonic if $d\alpha=0=\delta \alpha$, where
$\delta=(-1)^{k+1}*d*$ and $*$ is the symplectic star operator, see
Section 2. We set
$$
H^k_{\hr}(M)=H^k_{\hr}(M,\omega):=\Omega^k_{\hr}(M)/(\Im\, d\cap
\Omega^k_{\hr}(M))
$$
and
$$
h_k(M)=h_k(M,\omega):=\dim\,H^k_{\hr}(M,\omega).
$$
Since every symplectically harmonic form is closed, $H^k_{\hr}(M)$ is
a subgroup of $H^k(M)$ and $h_k\leq b_k$.

Of course, the above definition of the numbers $h_k$ is not symmetric:
one can also consider the ``dual'' numbers
$$
h^*_k(M)=h^*_k(M,\omega):=\dim\,\left(\Omega^k_{\hr}/\Im\,\delta
\cap\Omega^k_{\hr}\right).
$$

It turns out that the duality isomorphism $*:\Omega^{m-k}_{\hr} \to
\Omega^{m+k}_{\hr}$ yields the equality $h_{m-k}=h^*_{m+k}$, see Section 7.

Relating with the study of symplectically harmonic forms, we are
interested in the following question raised by Boris Khesin and Dusa
McDuff, see Yan~\cite{Yan}.

{\bf Question:} Which compact manifolds $M$ possess a continuous
family $\omega_t$ of symplectic forms such that $h_k(M,
\omega_t)$ varies with respect to $t$?

This question, according to Khesin, is probably related to group
theoretical hydrodynamics and geometry of diffeomorphism groups. Some
indirect indications for an existence of such relations can be found
in \cite{AK}.

Since we want to consider manifolds as in the question, it makes sense
to give them a certain name. So, let us call a closed smooth manifold
$M$ {\it flexible} if $M$ possesses a continuous family of symplectic
forms $\omega_t,\,t\in[a,b]$, such that $h_k(M,
\omega_a)\neq h_k(M, \omega_b)$ for some $k$. So, the above
question asks about the existence of flexible manifolds.

Yan \cite{Yan} has studied the case of closed $4$-manifolds. He proved
that 4-dimensional nilmanifolds are not flexible. He has also
found examples of flexible 4-manifolds. Actually, one step in this
proof of the existence is wrong, but the whole proof can easily be repaired, and hence the existence result holds. See Section 4 for details.

So, passing to higher dimensions, we have the following
question:

\noindent{\bf Question:} Do there exist flexible nilmanifolds of dimension
$\geq 6$?

Related with this question, we show that some $6$-dimensional
nilmanifolds are flexible.

Some words about tools. First, in order to prove the flexibility, we
must be able to compute the symplectically harmonic Betti numbers. It
turns out that, for every closed symplectic manifold $(M^{2m},\omega)$
and $k=1$ or 2 we have
$$
h_{2m-k}(M,\omega)=\rk\, (L_{[\omega]}^{m-k}: H^k(M) \to H^{2m-k}(M)),
$$
see \ref{yan} and \ref{tool}. So, a purely cohomological information
is enough in order to compute $h_{2m-k}$. Furthermore, if
$h_{2m-k}(M,\omega_1)\neq h_{2m-k}(M,\omega_2),
\, k=1,2$, then $M$ is flexible, see \ref{flex1} and \ref{yan}.

 In general, not only symplectically harmonic forms but also Poisson
harmonic forms are being of great interest in different areas of
mathematics and physics, \cite{EL}.

We use the sign \QQED\ in order to indicate the end of a proof. However, if we formulate a claim without proof, we put the sign $\square$ in the end of the claim.

\section{Symplectically harmonic forms}
\setcounter{equation}{0}

Let $(M^{2m},\omega)$ be a symplectic manifold. It is well known that
there exists a unique non-degenerate Poisson structure $\Pi$
associated with the symplectic structure (see, for example,
\cite{L,TO}), that is, $\Pi$ is a skew-symmetric tensor field of
order 2 such that $[\Pi,\Pi]=0$, where $[-,-]$ is the
Schouten-Nijenhuis bracket.

The Koszul differential $\delta:\Omega^k(M)\longrightarrow
\Omega^{k-1}(M)$ is defined for symplectic manifolds, and more
generally,
for Poisson manifolds as
$$
\delta=[i(\Pi),d].
$$

Brylinski have proved in \cite{Br} that the Koszul differential is a
symplectic codifferential of the exterior differential, with respect
to the symplectic star operator. We choose the volume form associated
to the symplectic form, that is, $v_M=\omega^m/m!$. Then we define the
symplectic star operator
$$
\ast:\Omega^k(M)\longrightarrow \Omega^{2m-k}(M),
$$
by the condition $\beta\wedge
(\ast\alpha)=\Lambda^k(\Pi)(\beta,\alpha)v_M$, for all
$\alpha,\beta\in \Omega^k(M)$. The symplectic star operator satisfies
the identities
$$
\ast^2=id,\quad \delta =(-1)^{k+1} \ast d\ast, \quad
\hbox{and}\quad i(\Pi)= L^{\ast}:=-\ast L\ast.
$$
Furthermore, if $M$ is a symplectic manifold $M$ then the operators $L$,
$d$, $\delta$ and $L^{\ast}$ (acting on the algebra $\Omega^{\ast}(M)$)
satisfy the following commutator relations:
\begin{equation}\label{relation}
[L,d]=0,\quad [L,\delta]=-d,\quad [L^{\ast},\delta]=0, \quad
[L^*,d]=-\delta.
\end{equation}

\begin{remark}\label{liberman}
{\rm The symplectic star operator was first considered by Libermann
\cite{L}
as $\ast (\alpha)=i(\mu^{-1}\alpha)v_M$, where $\mu$ is the
canonical isomorphism between the exterior algebras of vector fields and
forms. She has also introduced and studied the operators $\delta=(-1)^{k+1} \ast
d\ast$ and
$L^{\ast}=-\ast L\ast$. In particular, she has proved (using the Lepage
decomposition) that
\begin{equation}\label{LL}
\{ \alpha\in \Omega^{m-k}(M)\, |\, L^ {k+1}\alpha=0\}=
\{ \alpha\in \Omega^{m-k}(M)\, |\, L^ {\ast}\alpha=0\}.
\end{equation}
(See \cite{Yan} for a proof using the theory of
${\frak s\frak l}(2)$-representations.) }
\end{remark}

\begin{definition}{\rm
A $k$-form $\alpha$ on the symplectic manifold $M$ is called {\it
symplectically harmonic} if $d\alpha=\delta\alpha=0$. We denote
the space of harmonic $k$-forms by $\Omega^k_{\hr}(M)$. }
\end{definition}

It is clear from (\ref{relation}) that the form $\omega \wedge\alpha$ is
symplectically harmonic whenever $\alpha$ is. Hence, for every $k$ we
have the mapping $L:\Omega^k_{\hr}(M)\longrightarrow
\Omega^{k+2}_{\hr}(M)$.

We set
$$
H^k_{\hr}(M)=\Omega^{k}_{\hr}(M)/\Im\, d \cap \Omega^{k}_{\hr}(M)
{\rm\ and\ }h_k=h_k(M,\omega)=\dim\,H^k_{\hr}(M).
$$
So, for every symplectic manifold $M$, its de Rham cohomology
$H^{\ast}(M)$ contains a symplectically harmonic subspace
$H^{\ast}_{\hr}(M)$. We say that a de Rham cohomology class is {\it
symplectically harmonic} if it contains a symplectically harmonic
representative, i.e. if it belongs to the image of the inclusion
$H^*_{\hr}(M)
\subset H^*(M)$. Finally, we say that a manifold $M$ is {\it flexible} if
$M$ possesses a continuous family of symplectic forms
$\omega_t,\,t\in[a,b]$, such that $h_k(M,
\omega_a)\neq
h_k(M, \omega_b)$ for some $k$.

\begin{remark}
{\rm In the general case of a (degenerate) Poisson manifold
$(M,\Pi)$, we say that a $k$-form $\alpha$ is {\it Poisson
harmonic} if $d\alpha=0=\delta\alpha$. Notice that for a Poisson
manifold, in particular for a symplectic manifold,
$\Delta=d\delta+\delta d \equiv 0$, contrarily to the Riemannian case.
}
\end{remark}

\begin{remark}
{\rm Recall that a {\it symplectomorphism} between two symplectic manifolds
$(M,\omega_1)$ and $(N,\omega_2)$ is a diffeomorphism
$\phi:M\longrightarrow N$ such that $\phi^{\#}(\omega_2)=\omega_1$.
It is easy to see that
$$
\phi^{\ast} ( H^{k}_{\hr}(N))= H^{k}_{\hr}(M).
$$
where $\phi^{\ast}: H^{k}(N) \to H^{k}(M)$ is the induced homomorphism
in the de Rham cohomology. In other words, $H^*_{\hr}(-)$ is a
symplectic invariant. In particular, if $h_k(M,\omega_1)
\neq h_k(M,\omega_2)$ for two symplectic forms $\omega_1,\, \omega_2$
on $M$ then $\omega_1$ and $\omega_2$ are not symplectomorphic. }
\end{remark}

We do not know whether a smooth map (not a diffeomorphism) $\psi$ with
$\psi^{\#}\omega_2=\omega_1$ induces a map of symplectically harmonic
cohomology.

\begin{proposition}\label{dual}
{\rm (\cite{Yan})} For every
symplectic manifold $(M,\omega)$, the homomorphism
$$
L^k: \Omega^{m-k}_{\hr}(M)\longrightarrow \Omega^{m+k}_{\hr}(M)
$$
is an isomorphism.
\hqed
\end{proposition}

\begin{corollary}\label{sur} 
The homomorphism
$$
L^k: H^{m-k}_{\hr}(M)\longrightarrow H^{m+k}_{\hr}(M)
$$
is an epimorphism. In particular, $h_{m-k}\geq h_{m+k}$.
\hqed
\end{corollary}

\begin{corollary}\label{l:diag} Let $(M^{2m},\omega)$ be a symplectic
manifold. Then
$$
H^{m+k}_{\hr}(M)=\Im\{L^k: H^{m-k}_{\hr}(M) \to  H^{m+k}(M)\}
\subset H^{m+k}(M).
$$
\end{corollary}

\begin{proof}
It is a direct consequence of the commutativity of the  diagram
\begin{equation}
\diagram{ \Omega^{m-k}_{\hr}(M) & \hfl{L^k}{} &
\Omega^{m+k}_{\hr}(M)  \cr
\vvfl{}{} && \vvfl{}{} \cr
H^{m-k}_{\hr}(M) & \hfl{L^k}{} & H^{m+k}_{\hr}(M) \cr }
\end{equation}
since the top map $L^k$ is an isomorphism by \ref{dual} and both
vertical maps are the epimorphisms.
\end{proof}

\begin{corollary}
Let $(M^{2m},\omega)$ be a closed symplectic manifold. If
$h_{m+k}(M)=b_{m+k}(M)$ then $h_{m-k}(M)=b_{m-k}(M)$.
\hqed
\end{corollary}

\begin{proof}
Because of \ref{sur} and Poincar\'e duality, $h_{m-k}\geq h_{m+k}=b_{m+k}=b_{m-k}\geq h_{m-k}$.
\end{proof}

\begin{corollary}
Let $(M, \omega)$ be a symplectic manifold. If $b_k(M)=0$ for some
$k\leq m$, then $h_{2m-i}(M,\omega)=0$ for $i\leq k$ with $k-i$ even.
\end{corollary}

\begin{proof} It follows from \ref{l:diag}, since the homomorphism
$$
L^{m-i}:  H^i(M) \longrightarrow H^{2m-i}(M)
$$
passes through the trivial group $H^k(M)$.
\end{proof}

\begin{lemma}\label{linalg}
Let $\L$ be the space of all linear maps $\Bbb R^k \to \Bbb R^l$. Then
the following holds:
\par {\rm (i)} for every $r$ the set
$$
\{A\in \L\bigm | \rk\,A \leq r\}
$$
is an algebraic subset of $\L$. $($Here we regard $\L$ as the space
$\R^{kl}$ of $l\times k$-matrices whose entries are regarded as the
coordinates$)$;
\par {\rm (ii)} for every $m\leq \min\{k,l\}$ the set
$
\{A\in \L\bigm | \rk\,A \geq m\}
$
is open and dense in $\L$;
\par {\rm(iii)} let $A,B\in \L$ be two linear maps such that $\rk\, A <\rk\, B$. Then the set
$$
\Lambda=\{\lambda \in \R \bigm | \rk\,(A+\lambda B)\geq \rk\, B\}
$$
is an open and dense subset of $\R$.
\end{lemma}

\begin{proof}
 (i) This claim follows, because the rank of a matrix is equal to the
 order of the largest non-zero minor.
 \par (ii) This claim follows from (i).
 \par (iii) By (i), the set $\R \setminus \Lambda$ is an algebraic
 subset of $\R$. So, it suffices to prove that $\Lambda\neq
 \emptyset$. But, by (ii), $\rk\, (B+\mu A)> \rk\, A $ for $\mu$ small enough, and
 so $\Lambda\neq\emptyset$.
\end{proof}

\begin{corollary}\label{flex1}
Let $\omega_0$ and $\omega_1$ be two symplectic forms on a manifold
$M^{2m}$. Suppose that, for some $k>0$, $h_{m-k}(\omega_0)=
h_{m-k}(\omega_1)$, but $h_{m+k}(\omega_0)< h_{m+k}(\omega_1)$. Then,
for every $\eps>0$, there exists $\lambda\in (0,\eps)$ such that
$h_{m+k}(\omega_0+\lambda
\omega_1)> h_{m+k}(\omega_0)$. Moreover, $M$ is flexible
provided that it is closed.
\end{corollary}

\begin{proof}
The existence of $\lambda$ follows from \ref{l:diag} and
\ref{linalg}(iii). The flexibility of $M$ follows, since
$\omega_0+t\omega_1$ is a symplectic form for $t$ small enough. Indeed, if
we set $\omega_t=\omega_0+t\omega_1, t\in [0,\lambda]$, then
$h_{m+k}(\omega_0) < h_{m+k}(\omega_{\lambda})$.
\end{proof}

\begin{corollary}\label{flex2}
Let $(M^{2m}, \omega_0)$ be a closed
symplectic manifold. Given $k$ with $0<k<m$, suppose that
$h_{m-k}(M,\omega)= b_{m-k}(M)$ for every symplectic form $\omega$ on
$M$. Furthermore, suppose that there exists $x\in H^2(M)$ such that
$$
\rk\,\{L^k_x: H^{m-k}(M) \to H^{m+k}(M)\} > h_{m+k}(M,\omega_0).
$$
Then $M$ is flexible.
\end{corollary}

\begin{proof} Take a closed 2-form $\alpha$ which represents $x$. Then
$\omega_0+t\alpha$ is a symplectic form for $t$ small enough.
Furthermore, by \ref{l:diag} and \ref{linalg}(iii), there exists
arbitrary small $\lambda$ such that $h_{m+k}(\omega_0+\lambda
\alpha)> h_{m+k}(\omega_0)$. Now the result follows from
\ref{flex1}.
\end{proof}

We set
$$
\osi (M)=\{\omega\in \Omega^2(M)\bigm|\omega\ \hbox{is a symplectic form on}\ M\}
$$
and define $\Omega(b,k)=\{\omega\in \osi\bigm|h_k(M,\omega)=b\}$.

\begin{corollary}\label{strat}
Let $M^{2m}$ be a manifold that admits a symplectic structure. Suppose that, for some $k>0$, $h_{m-k}(M,\omega)$
does not depend on the symplectic structure $\omega$ on $M$. Then the following
three conditions are equivalent:
\par {\rm (i)} the set $\Omega (b, m+k)$ is open and dense in $\osi(M)$;
\par {\rm (ii)} the interior of the set $\Omega (b, m+k)$ in $\osi(M)$ is non-empty;
\par {\rm (iii)} the set $\Omega (b, m+k)$ is non-empty and
$h_{m+k}(M,\omega)\leq b$ for every $\omega\in \osi(M)$.
\end{corollary}
\begin{proof} (i) $\Rightarrow$ (ii). Trivial.
\par (ii) $\Rightarrow$ (iii). Suppose that there exists $\omega_0$ with
$h_{m+k}(M, \omega_0)> b$. Take $\omega$ in the interior of
$\Omega(b,m+k)$. Then, in view of \ref{linalg}, there exists an
arbitrary small $\lambda$ such that $h_{m+k}(\omega+\lambda
\omega_0)>b$, i.e. $\omega$ does not belong to the interior of
$\Omega (b,m+k)$. This is a contradiction.
\par (iii) $\Rightarrow$ (i). This is true because of \ref{linalg}(ii).
\end{proof}

So, the family $\{\Omega(b,m+k)|b=0,1, \ldots \}$ gives us a
stratification of $\osi (M)$ where the maximal strat is open and
dense.

\begin{lemma}\label{rho}
Let $(M^{2m},\omega)$ be a closed symplectic manifold, and set
$$
\rho_{2k+1}=\rk \left\{L^{m-2k-1}:H^{2k+1}(M)\longrightarrow
H^{2m-2k-1}(M)\right\}, \quad k=0,\ldots ,\left[{m-1\over 2}\right].
$$
Then $\rho_{2k+1}$ is an even number. Furthermore, $h_{2m-2k-1}\leq
\rho_{2k+1}\leq b_{2k+1}$, and $\rho_{2k+1} = b_{2k+1}$ if and only if
$L^{m-2k-1}:  H^{2k+1}(M) \longrightarrow H^{2m-2k-1}(M) $ is
surjective.
\end{lemma}

\begin{proof}
Let $p:H^{2k+1}(M)\otimes H^{2m-2k-1}(M)\longrightarrow \R$ be the usual
non-singular pairing given by
$$
p\left([\alpha],[\gamma]\right)=\int_M \alpha\wedge
\gamma,
$$
for $[\alpha]\in H^{2k+1}(M)$ and $[\gamma]\in
H^{2m-2k-1}(M)$. Define a skew-symmetric bilinear form $\langle -,
-\rangle :H^{2k+1}(M)\otimes
H^{2k+1}(M)\longrightarrow \R$ via the formula
$$
 \left\langle [\alpha],[\beta] \right\rangle
=p\left([\alpha],L^{m-2k-1}[\beta]\right),
$$
for $[\alpha],[\beta]\in H^{2k+1}(M)$. It is easy to see that the
rank of
$\langle -,- \rangle $, which must be an even number $2l$ with
$0\leq 2l\leq b_{2k+1}$, is
equal   to $\rho_{2k+1}$.

The inequality $h_{2m -2k -1}\leq \rho_{2k+1}$ follows from
\ref{l:diag}.

The last claim holds since, by the Poincar\'e duality,
$b_{2k+1}(M)=b_{2m-2k-1}(M)$.
\end{proof}

\begin{remark}\label{even}
{\rm Lemma \ref{rho} yields the following well-known fact: if $M$ in
\ref{rho} satisfies the Hard Lefschetz condition, then all
odd-dimensional Betti numbers are even. }
\end{remark}

\begin{corollary}
Let $(M^{2m},\omega)$ be a closed symplectic manifold. If
$b_{2k+1}(M)$ is odd then $h_{2m-2k-1}<b_{2m-2k-1}=b_{2k+1}$. In
particular, if $b_{2k+1}=1$ then $h_{2m-2k-1}=0$.
\hqed
\end{corollary}

\begin{corollary}\label{indep}
If $b_{4k+2}(M)=1$ then $\rho_{2k+1}$ does not depend on symplectic
structure on $M$.
\end{corollary}

\begin{proof} Because of the Poincar\'e duality, $b_{2m-4k-2}=1$. So, $[\omega]^{m-2k-1}$, and hence
$$
L^{m-2k-1}: H^{2k+1}(M) \longrightarrow H^{2m-2k-1}(M),
$$
is determined uniquely up to non-zero multiplicative constant.
\end{proof}

\section{The numbers $h_k$ and $h_{2m-k}$ for $k$ small}
\setcounter{equation}{0}

\begin{proposition}\label{harm}
Let $(M,\omega)$ be a symplectic manifold, and let $k$ be a non-negative integer
number such that the following holds:
\par {\rm (i)} $L^{k+2}:H^{m-k-2}(M)\longrightarrow H^{m+k+2}(M)$
is surjective;
\par {\rm (ii)} If $a\in H^{m-k-2}(M)$ is not symplectically
harmonic, then $La=0$.

Then every cohomology class in $H^{m-k}(M)$ is symplectically harmonic.
\end{proposition}

\begin{proof} Here we use some ideas from \cite{Yan}. It follows from (i) that
$$
H^{m-k}(M)=\Im\, L + P_{m-k},
$$
where $P_{m-k}=\{a\in H^{m-k}(M) \, |\, L^{k+1}a=0\}$. Indeed, if
$a\in H^{m-k}(M)$ then there exists $b\in H^{m-k-2}(M)$ such that
$L^{k+1}a=L^{k+2}b$. Therefore,
$a-b\wedge [\omega]\in P_{m-k}$, and
$$a=b\wedge [\omega] + (a-b\wedge [\omega])\in \Im\, L + P_{m-k}.$$

Because of (ii), every class in $\Im\, L$ is symplectically harmonic.
So, it suffices to prove that any cohomology class $a\in P_{m-k}$ is
symplectically harmonic.

Let $a=[\alpha]\in P_{m-k}$ with $\alpha\in \Omega^{m-k}(M)$
closed. Since $L^{k+1}a=0\in H^{m+k+2}(M)$, there exists $\beta\in
\Omega^{m+k+1}(M)$ such that $\alpha\wedge \omega^{k+1}= d\beta$. It
is known \cite{L} that $L^{k+1}:\Omega^{m-k-1}(M)\longrightarrow
\Omega^{m+k+1}(M)$ is
surjective. So, there exists $\gamma\in \Omega^{m-k-1}(M)$ with
$\beta=\gamma\wedge \omega^{k+1}$, and hence $(\alpha-d\gamma)\wedge
\omega^{k+1}=0$. So,  if we take $\overline{\alpha}=\alpha-d\gamma$,
then $[\overline{\alpha}]=a$ and $L^{k+1}\overline{\alpha}=0$. But the
equality $L^{k+1}\overline{\alpha}=0$ implies that
$L^{\ast}\overline{\alpha}=0$ in view of (\ref{LL}). Thus,
$\overline{\alpha}$ is symplectically harmonic by (\ref{relation}).
\end{proof}

\begin{corollary}\label{yan}
{\rm (\cite{Yan})} Let $(M,\omega)$ be an arbitrary symplectic manifold.
Then every cohomology class in $H^{k}(M),\, k=0,1,2,$ is
symplectically harmonic.
\hqed
\end{corollary}

Recall that a symplectic manifold $(M^{2m},\omega)$ is called a
{\it manifold of Lefschetz type}, if the map
$$
L:H^1(M)\longrightarrow H^{2m-1}(M)
$$
is surjective. Notice that, similarly to \ref{even}, $b_1(M)$ is even
for every closed manifold $M$ of Lefschetz type.

\begin{corollary}\label{lefschetz}
Let $(M^{2m},\omega)$ be a manifold of Lefschetz type. Then
every cohomology class in $H^3(M)$ is symplectically harmonic.
\hqed
\end{corollary}

\begin{corollary}\label{tool}
For every symplectic manifold $(M^{2m},\omega)$ and $k=0,1,2$,
$$
H^{2m-k}_{\hr}(M)=\Im\{L^{m-k}: H^{k}(M) \to H^{2m-k}(M)\} \subset
H^{2m-k}(M).
$$
Furthermore, if $M$ is a Lefschetz type manifold then
$$
H^{2m-3}_{\hr}(M)=\Im\{L^{m-3}: H^{3}(M)\to  H^{2m-3}(M)\} \subset
H^{2m-3}(M).
$$
\end{corollary}

\begin{proof} This follows from \ref{l:diag} because of \ref{yan}
and \ref{lefschetz}.
\end{proof}

\begin{corollary}
For every closed symplectic manifold $(M^{2m},\omega)$ the number
$h_{2m-1}(M^{2m},\omega)$ is even. Furthermore, if $b_2=1$ then
$h_{2m-1}(M^{2m},\omega)$ does not depend on $\omega$.
\end{corollary}

\begin{proof}
It follows from \ref{rho} and \ref{indep} since, by \ref{tool},
$\rho_1=h_{2m-1}$.
\end{proof}

\begin{corollary}
If $(M^{2m},\omega)$ is a symplectic manifold that is not a manifold of Lefschetz
type, then $h_{2m-1}\leq b_1-1$. In particular, if $b_1=2$, then $
h_{2m-1}=0$.
\hqed
\end{corollary}

This is the case for compact non-toral nilmanifolds (see \cite{TO} and
the table of the classification of $6$-dimensional compact
nilmanifolds in Section 5).

\section{Yan's result on flexibility in dimension 4}
\setcounter{equation}{0}

According to Yan \cite{Yan}, $4$-dimensional compact nilmanifolds are
not flexible. Indeed, by \ref{yan}, only $h_3$ may vary, but it turns
out that $h_3$ is constant. Namely, based on certain results from
\cite{FGC}, Yan~\cite{Yan} noticed the following relation for closed
4-dimensional nilmanifolds:
\begin{enumerate}
\item if $b_1(M)=2$, then $h_3=0$,
\item if $b_1(M)=3$ (therefore, $b_2(M)=4$),  then $h_3=2$,
\item if $b_1(M)=4$ (i.e. $M=\T^4$), then $h_3=4$.
\end{enumerate}

On the other hand, Yan~\cite{Yan} has found closed 4-dimensional
flexible manifolds, although his arguments need a certain correction
(see below). Namely, he formulated without proof the following
proposition, where $M$ is assumed to be a closed 4-dimensional manifold.

\begin{proposition}\label{4.1}
{\rm (\cite[Prop. 4.1]{Yan})}
The following assertions are equivalent:
\par {\rm (i)} There exists a family $\omega_t$ of symplectic forms such that
$h_3$ varies.
\par {\rm (ii)}  There exist two symplectic forms $\omega_1$ and $\omega_2$
such that $\Im\, L_{[\omega_1]}\not= \Im\, L_{[\omega_2]}$, where
$L_{[\omega_i]}$ is the Lefschetz map respect to $[\omega_i]$,
$(i=1,2)$.
\par {\rm (iii)} There exists a symplectic form $\omega$ on $M$ and a class
$a\in H^2(M)$ such that $\Im\, L_{a}\not\subset \Im\, L_{[\omega]}$.
\par {\rm (iv)} There exists a symplectic form $\omega$ on $M$ such that $\Im\,
L_{[\omega]}$ is not equal to the image of the cup product pairing
$H^1(M)\otimes H^2(M)\longrightarrow H^3(M)$.
\hqed
\end{proposition}

Concerning to this proposition, it is true that (i)$\Rightarrow$ (ii)
$\Rightarrow$ (iii) $\Rightarrow$ (iv), but the Kodaira--Thurston
manifold satisfies (iv) and does not satisfy (i). The Kodaira--Thurston
manifold is obtained by taking the product of the Heisenberg manifold
and the circle (this manifold is a compact nilmanifold). Its Sullivan
minimal model has the form
$$
(\Lambda(x_1,x_2,x_3,x_4),d)\quad \hbox{with}\quad dx_1=dx_2=dx_4=0,\,
dx_3=x_1x_2,
$$
where $\deg\,x_i=1$ for $i=1,2,3,4$ and the generators $x_1, x_2, x_3$
come from the Heisenberg manifold. The cohomology class of the
symplectic form is then given by the element $\omega=x_1x_4+x_2x_3$.
Now, from a direct computation we obtain that $\Im\, L_{[\omega]}$ is
generated by the non-zero classes of $x_1x_2x_3$ and $x_2x_3x_4$ and
the image of the cup product $H^1(M)\otimes H^2(M)\longrightarrow
H^3(M)$ is generated by $x_1x_2x_3$, $x_1x_3x_4$ and $x_2x_3x_4$, so
condition (iv) is satisfied. Moreover, it is easy to see that the
Kodaira--Thurston manifold satisfies the condition (ii).

But, because of what we have said in the beginning of the section, any 4-dimensional nilmanifold (and hence the Kodaira-Thurston
manifold) does not satisfy condition (i).

Yan's construction of flexible closed 4-dimensional manifolds is
based on the following proposition.

\begin{proposition} \label{4-flex}
{\rm \cite[Cor. 4.2]{Yan}} Let $(M^4,\omega)$ be a closed symplectic
manifold which satisfies the following conditions:
\par {\rm (i)} the homomorphism $L_{[\omega]}:H^1(M)\longrightarrow H^3(M)$ is
trivial;
\par {\rm (ii)} the cup product $H^1(M)\otimes H^2(M)\longrightarrow H^3(M)$
is non-trivial.

Then $M$ is flexible.
\hqed
\end{proposition}

Yan regards this proposition as a corollary of his Proposition 4.1. As
we have seen, the last one is wrong. However, Proposition \ref{4-flex}
is correct because it is a special case of our Corollary \ref{flex2}.

Finally, Gompf \cite[Observation 7]{G} proved the existence of
4-manifolds as in \ref{4-flex}.

\section{Symplectically harmonic forms in homogeneous $\quad$ spaces}
\setcounter{equation}{0}

From the previous section we know that none 4-dimensional nilmanifold
is flexible. The goal of this section is to demonstrate the existence
of 6-dimensional flexible nilmanifolds.

A compact nilmanifold is a homogeneous space of the form $G/\Gamma$,
where $G$ is a simply connected nilpotent Lie group and $\Gamma$ is a
discrete co-compact subgroup of $G$, i.e. a lattice (co-compact means
that $G/\Gamma$ is compact). Recall that $\Gamma$ is determined by $G$
uniquely up to an isomorphism. In greater detail, if $\Gamma$ and
$\Gamma'$ are two lattices in $G$ then there exists an automorphism
$\varphi: G \to G$ with $\varphi(\Gamma)=\Gamma'$, see~\cite{R}.
Moreover, $\Gamma$ determines $G$ uniquely up to an isomorphism. In
particular, the compact nilmanifold $G/\Gamma$ determines and is
completely determined by $G$.

Three important facts in the study of compact nilmanifolds are (see
\cite{TO}):
\begin{enumerate} \item Let $\g$ be a nilpotent Lie algebra with
structural constants $c^{ij}_k$ with respect to some basis, and let $\{\alpha_1,\ldots
,\alpha_n\}$ be the dual basis of $\g^{\ast}$. Then in the Chevalley--Eilenberg complex
$(\Lambda^{\ast}\g^{\ast},d)$ we have
\begin{equation}\label{cte}
d\alpha_k=\sum_{1\leq i<j<k} c^{ij}_k
\alpha_i\wedge \alpha_j.
\end{equation}
\item Let $\g$ be the Lie algebra of a simply connected nilpotent Lie group $G$. Then, by Malcev's theorem, $G$ admits a lattice if
and only if $\g$ admits a basis such that all the structural constants
are rational.
\item By Nomizu's theorem, the Chevalley--Eilenberg complex
$(\Lambda^{\ast}\g^{\ast},d)$ of $\g$ is quasi-isomorphic to the de Rham
complex of $G/\Gamma$. In particular,
\begin{equation}\label{nom}
H^{\ast}(G/\Gamma)\cong H^{\ast}(\Lambda^{\ast}\g^*,d)
\end{equation}
and any cohomology class $[a]\in H^k(G/\Gamma)$
contains a homogeneous representative $\alpha$. Here we call the form $\alpha$
homogeneous if the pullback of $\alpha$ to $G$ is left invariant.
\end{enumerate}

\begin{theorem}\label{nil}
There exist at least five $6$-dimensional flexible nilmanifolds.
\end{theorem}

To prove this theorem we run our fingers over all 34 $6$-dimensional
compact nilmanifolds. The results are contained in the table below.
The proofs take all the remained part of the section.

It follows from Corollary \ref{yan} that $\dim\, H^k_{\hr}(M)=\dim\,
H^k(M)$, for $k\not=3,4,5$. Therefore, we should study the behaviour
of $h_k=\dim\, H^k_{\hr}(M)$ for degrees $k=3,4,5$. We were not able
to compute $h_3$, but we have found 5 manifolds with $h_4$ and/or
$h_5$ varying. So, by \ref{flex1}, there are at least five
6-dimensional flexible nilmanifolds. (There are some reasons to
conjecture that $h_3=b_3$ for all closed 6-dimensional manifolds. So,
if it is true then we have exactly five 6-dimensional flexible
nilmanifolds.)

It turns out that every 6-dimensional real nilpotent Lie algebra
admits  a basis with rational structural constants.  So, by what we
said above, the 6-dimensional compact nilmanifolds  are in a bijective
correspondence with the 6-dimensional simply connected nilpotent Lie
groups, and hence with the 6-dimensional nilpotent Lie algebras.

We use the classification of nilpotent Lie algebras given by Salamon~
\cite{Sa}. It is based on the Morozov classification of
$6$-dimensional nilpotent Lie algebras \cite{OV}. We have added to
Salamon's classification the symplectically harmonic Betti numbers $
h_k(M)$ for $k=4,5$.

In the table Lie algebras appear lexicographically with respect to
$(b_1,b_2,6-s)$ where $s$ is the step length. The first two columns
contain the Betti numbers $b_1$ and $b_2$ (notice that
$b_3=2(b_2-b_1+1)$ because of the vanishing of the Euler
characteristic). The next column contains $6-s$, where $s$ is the step
length.
\par The fourth column contains the description of the structure of the
Lie algebra by means of the expressions of the form (\ref{cte}) in the
Chevalley-Eilenberg complex. In view of \ref{nom}, it means that, say,
for the compact nilmanifold $M$ from the second row, there exists a
basis $\{\alpha_i\}_{i=1}^6$ of homogeneous $1$-forms on $M$ such that
$$
d \alpha_1=0=d \alpha_2,\ d \alpha_3=\alpha_1 \wedge \alpha_2,\ d
\alpha_4=\alpha_1 \wedge \alpha_3,\ d \alpha_5=\alpha_1 \wedge
\alpha_4,\ d \alpha_6=\alpha_3 \wedge
\alpha_4 + \alpha_5 \wedge \alpha_2.
$$
\par The column headed $\oplus$ indicates the dimensions of the
irreducible subalgebras in case $\g$ is not itself irreducible.
\par The next columns show the dimensions
$h_k$ for $k=4,5$. So, the column, say, $h_4$ contains all possible
values of $h_4(M,\omega)$ which appear when $\omega$ runs over all
symplectic forms on $M$. The sign ``--'' at a certain row means that
the corresponding Lie algebra (as well as the compact nilmanifold)
does not admit a symplectic structure.
\par For completeness, in the last columns we list the real
dimension $\dim_{\R} \SS(\g)$ of the moduli space of symplectic structures.

\begin{convention}{\rm
\begin{enumerate} \item
From now on we  write $\alpha_{ij\cdots\,k}$
instead of $\alpha_i\wedge \alpha_j\wedge\cdots\wedge
\alpha_k$.
\item
In future we say that a compact nilmanifold $G/\Gamma$ has type, say
(0,0,12,13,14,15) if the corresponding Lie algebra has the
structure (0,0,12,13,14,15) (i.e. in our case, sits in the
third row).
\end{enumerate}
}
\end{convention}

\vfil\eject\begin{center}
{\large\bf Six-dimensional real nilpotent Lie algebras}
%\vspace{20pt}

\def\no{\hb{\bf--}}\[\ba{|c|c|c|l|c|c|c|c|}\hline
b_1&b_2&6\!-\!s&\hb{\qq Structure}&\op& h_4& h_5&\dim_{\R}
\SS(\g)\\\hline 2&2&1&(0,0,12,13,14+23,34+52)&  &\no&\no& \no\\
2&2&1&(0,0,12,13,14,34+52)&  &\no&\no&\no\\ 
2&3&1&(0,0,12,13,14,15)&
&3&0& 7\\ 2&3&1&(0,0,12,13,14+23,24+15)&  &2&0& 7\\
2&3&1&(0,0,12,13,14,23+15)&  &2&0& 7\\ 
2&4&2&(0,0,12,13,23,14)&  &4&0&
8\\ 2&4&2&(0,0,12,13,23,14-25)&  &2,3,4&0& 8\\
2&4&2&(0,0,12,13,23,14+25)&  &4&0& 8\\\hline
3&4&2&(0,0,0,12,14-23,15+34)&  &2&0& 7\\ 3&5&2&(0,0,0,12,14,15+23)&
&4&2& 8\\ 3&5&2&(0,0,0,12,14,15+23+24)&  &3,4&0,2& 8\\
3&5&2&(0,0,0,12,14,15+24)&1+5&4 &2 & 8\\ 3&5&2&(0,0,0,12,14,15)&1+5&4
&2 & 8\\ 3&5&3&(0,0,0,12,13,14+35)&  &\no&\no&\no\\
3&5&3&(0,0,0,12,23,14+35)&  &\no&\no& \no\\ 3&5&3&(0,0,0,12,23,14-35)&
&\no &\no &\no\\ 3&5&3&(0,0,0,12,14,24)&1+5&\no &\no &\no\\
3&5&3&(0,0,0,12,13+42,14+23)&  &3 &0 & 8\\ 3&5&3&(0,0,0,12,14,13+42)&
&3 &0 & 8\\ 3&5&3&(0,0,0,12,13+14,24)&  &2,3 &0 &8\\
3&6&3&(0,0,0,12,13,14+23)&  &3 &0 & 9\\ 3&6&3&(0,0,0,12,13,24)&  &5 &0
& 9\\ 3&6&3&(0,0,0,12,13,14)&  &4&0 & 9\\ 3&8&4&(0,0,0,12,13,23)&
&7,8&0 & 9\\\hline 4&6&3&(0,0,0,0,12,15+34)&  &\no&\no&\no\\
4&7&3&(0,0,0,0,12,15)&1\!+\!1\!+\!4&3&2& 9\\
4&7&3&(0,0,0,0,12,14+25)&1+5&3&2& 9\\ 4&8&4&(0,0,0,0,13+42,14+23)&
&7&2& 10\\ 4&8&4&(0,0,0,0,12,14+23)&  &6 &2 & 10\\
4&8&4&(0,0,0,0,12,34)&3+3&7 &2  & 10\\ 4&9&4&(0,0,0,0,12,13)&1+5&7,8
&2 & 11\\\hline 5&9&4&(0,0,0,0,0,12+34)&1+5&\no &\no & \no\\
5&11&4&(0,0,0,0,0,12)&1\!+\!1\!+\!1\!+\!3& 9 &4 & 12\\\hline
6&15&5&(0,0,0,0,0,0)&1+\cdots+1& 15 &6& 15\\\hline\ea\]
\end{center}\vfil\eject

{\bf Proof of Theorem \ref{nil}}

We prove the theorem via considering  case by case. Namely, we study
in more detail the cases which are proclaimed to be flexible. In view
of \ref{flex1}, they are precisely the cases of compact nilmanifolds with
varying symplectically harmonic Betti numbers $h_k$. Here the main
tool for computing $h_k$ is Corollary \ref{tool}. In order to find
symplectic structures on $M$, we use the following proposition.

\begin{proposition}\label{sympl}
Let $M^{2n}$ be a compact manifold of the form $G/\Gamma$ where
$\Gamma$ is a discrete subgroup of a Lie group $G$, and let $\omega\in
\Omega^2(M)$ be a closed homogeneous $2$-form such that
$[\omega]^n\neq 0$. Then $\omega$ is a symplectic form on $M$.
\end{proposition}

\begin{proof} Since $[\omega]^n\neq 0$, we conclude that the linear form
$\omega|T_xM$ is non-degenerate for some point $x\in M$. So, $\omega^n$ is
non-degenerate since it is homogeneous. Thus, $\omega$ is non-degenerate.
\end{proof}

\medskip

\begin{proposition}\label{case1}
The compact nilmanifold $M$ of the type $(0,0,12,13,23,14-25)$ is flexible.
\end{proposition}

\begin{proof}
According to our assumption about the type of $M$, there exists a basis
$\{\alpha_i\}_{i=1}^6$ of homogeneous $1$-forms on $M$ such that
$$
d\alpha_1=d\alpha_2=0,\quad d\alpha_3=\alpha_{12},\quad
d\alpha_4=\alpha_{13},\quad d\alpha_5=\alpha_{23},\quad
d\alpha_6=\alpha_{14}-\alpha_{25}.
$$

Since the de Rham cohomology of the nilmanifold is isomorphic to the
Chevalley-Eilenberg cohomology of the Lie algebra, we conclude that
\begin{eqnarray*}
H^1(M)&=&\{[\alpha_1],[\alpha_2]\},\\ H^2(M)&=&\{[\alpha_{14}],
[\alpha_{15}+\alpha_{24}], [\alpha_{26}-\alpha_{34}],
[\alpha_{16}-\alpha_{35}]\}.
\end{eqnarray*}

In particular, by \ref{sympl}, a 2-form $\omega$ on $M$ is symplectic if and
only if 
\begin{equation}
[\omega]=A[\alpha_{14}]+B [\alpha_{15}+\alpha_{24}]+
C [\alpha_{26}-\alpha_{34}]+D[\alpha_{16}-\alpha_{35}],
\end{equation}
where $ACD-B(C^2+D^2)\not=0$, $A,B,C,D\in \R$.

\begin{claim}\label{example1}
Let $\omega$ be a symplectic form on $M$. Then the following holds:
\par  {\rm (i)} if $C^2\not=D^2$, then $h_4=4$; 
\par {\rm
(ii)} if $C^2=D^2$ and $A^2\neq 4B^2$, then $h_4=3$;
\par{\rm(iii)} if $C^2=D^2$ and $A^2=4B^2$, then $h_4=2$.
\par
Furthermore, $h_5=0$ for every symplectic form $\omega$ on $M$.
\end{claim}

\begin{proof}
It follows from the Poincar\'e duality that $H^4(M)=\Bbb R^4$. Hence, by (\ref{nom})
$$
H^4(M)=\{[\alpha_{1246}], [\alpha_{1256}], [\alpha_{1356}],
[\alpha_{1346}+\alpha_{2356}]\}
$$
since the four cohomology classes from above are linearly independent.
Furthermore, for every $\omega$ the image of the mapping
$L:H^2(M)\longrightarrow H^4(M)$ is
\begin{eqnarray*}
\Im\, L &=& \{ -C[\alpha_{1246}]+D[\alpha_{1256}],
2D[\alpha_{1246}]-2C[\alpha_{1256}], \\
&&-A[\alpha_{1246}]-2B[\alpha_{1256}]-2C[\alpha_{1356}]-D[\alpha_{1346}
+\alpha_{2356}],\\ &&
2B[\alpha_{1246}]+A[\alpha_{1256}]-2D[\alpha_{1356}]
-C[\alpha_{1346}+\alpha_{2356}]\},
\end{eqnarray*}
which has dimension $4$ for $C^2\not= D^2$, dimension $3$ for $C^2=D^2$ and $A^2\neq 4B^2$, and dimension $2$ for $C^2= D^2$ and $A^2=4B^2$.
The result follows from Corollary \ref{tool}.
\end{proof}

Now, the proof of the proposition follows from \ref{flex1}.
\end{proof}

\begin{proposition}\label{case2}
The compact nilmanifold of the type $(0,0,0,12,14,15+23+24)$ is flexible.
\end{proposition}

\begin{proof} According to our assumption about the type of $M$, there exists
a basis $\{\alpha_i\}_{i=1}^6$ of homogeneous $1$-forms on $M$ such
that $$ d\alpha_1=d\alpha_2=d\alpha_3=0,\quad
d\alpha_4=\alpha_{12},\quad d\alpha_5=\alpha_{14},\quad
d\alpha_6=\alpha_{15}+\alpha_{23}+\alpha_{24}.
$$

The cohomology groups of degrees 1, 2 are:
\begin{eqnarray*}
H^1(M)&=&\{[\alpha_1],[\alpha_2], [\alpha_3]\},\\
H^2(M)&=&\{[\alpha_{13}], [\alpha_{15}], [\alpha_{23}],
[\alpha_{16}+\alpha_{25}-\alpha_{34}], [\alpha_{26}-\alpha_{45}]\}.
\end{eqnarray*}

In particular, if $\omega$ is a  symplectic form on
$M$ then
\begin{equation}
[\omega]=A [\alpha_{13}]+ B [\alpha_{15}]+ C
[\alpha_{23}]+D [\alpha_{16}+\alpha_{25}-\alpha_{34}]+E
[\alpha_{26}-\alpha_{45}],
\end{equation}
where $AE^2+BDE-CDE-D^3\not=0$.

The following claim can be proved similarly to \ref{example1}.

\begin{claim}\label{example2}
Let $\omega$ be a symplectic form on $M$. Then the following holds:
\par {\rm (i)} if $E\not= 0$, then $h_4=4$,
$h_5=2$;
\par
{\rm(ii)} if $E=0$, then $h_4=3$, $h_5=0$.
\hqed
\end{claim}

Now, by \ref{flex1}, $M$ is flexible.
\end{proof}

We hope that now it is clear how to run over all the three remaining cases.
So, below we omit the details while indicate the main steps of the corresponding
calculations.

\begin{proposition}\label{case345}
The compact nilmanifolds of types

$(0,0,0,12,13+14,24)$,
$(0,0,0,12,13,23)$ and $(0,0,0,0,12,13)$

are flexible.
\end{proposition}

\begin{proof}
\underline{\bf Case $(0,0,0,12,13+14,24)$:}
\begin{eqnarray}
H^2(M)=\{[\alpha_{13}], [\alpha_{15}], [\alpha_{23}], [\alpha_{16}+\alpha_{25}
+\alpha_{34}], [\alpha_{26}]\},
\end{eqnarray}
and the 2-form $\omega$ on $M$ is symplectic if and only if
$$
[\omega]=A[\alpha_{13}]+B [\alpha_{15}]+C [\alpha_{23}]+
D [\alpha_{16}+\alpha_{25}+\alpha_{34}]+E [\alpha_{26}],
$$
where $D(BE-D^2)\not= 0$. Furthermore, if $EB+3D^2\not= 0$ then $
h_4=3$; otherwise, $h_4=2$. Finally, $h_5=0$ for every $\omega$.

\medskip
\underline{\bf Case $(0,0,0,12,13,23)$:}
\begin{eqnarray*}
H^2(M)=\{[\alpha_{14}], [\alpha_{15}], [\alpha_{16}+\alpha_{25}],
 [\alpha_{16}-\alpha_{34}], [\alpha_{24}], [\alpha_{26}], [\alpha_{35}],
[\alpha_{36}]\},
\end{eqnarray*}
and the 2-form $\omega$ on $M$ is symplectic if and only if
$$
[\omega]=A[\alpha_{14}]+B [\alpha_{15}]+C [\alpha_{16}+\alpha_{25}]+D
 [\alpha_{16}-\alpha_{34}]+E [\alpha_{24}]+F [\alpha_{26}]+G [\alpha_{35}]+H
[\alpha_{36}],
$$
where
$$
ACH-AFG-BDF-BEH+DC^2+CEG+CD^2+DEG\not= 0.
$$
Furthermore, if $C^2+CD+D^2-BF-EG+AH\not=0$, then $h_4=8$; otherwise,
$h_4=7$. Finally, $h_5=0$ for every $\omega$.

\medskip
\underline{\bf Case $(0,0,0,0,12,13)$:}

The second de Rham cohomology group is given by
\begin{eqnarray*}
H^2(M)=\{[\alpha_{14}], [\alpha_{15}], [\alpha_{16}],[\alpha_{23}],
 [\alpha_{24}], [\alpha_{25}], [\alpha_{34}], [\alpha_{26}+\alpha_{35}],
[\alpha_{36}]\}.
\end{eqnarray*}
The 2-form $\omega$ on $M$ is symplectic if and only if
\begin{eqnarray*}
[\omega]=A[\alpha_{14}] &+& B [\alpha_{15}]+C
[\alpha_{16}]+D[\alpha_{23}]+ E [\alpha_{24}]+F [\alpha_{25}]\\
&+&G[\alpha_{34}]+H [\alpha_{26}+
\alpha_{35}]+ I[\alpha_{36}],
\end{eqnarray*}
where $-AFI+H^2A+BEI-BGH-CEH+CFG\not= 0$. Furthermore, if
$H^2-FI\not=0$, then $h_4=8$; otherwise, $h_4=7$. Finally, $h_5=2$ for
every $\omega$.

Notice that in this last case the nilpotent Lie algebra (so the
compact nilmanifold) is reducible of type $(1+5)$.
\end{proof}

\begin{comment}
{\rm Here we want to say more about flexibility of manifolds appeared
in \ref{case1} and \ref{case2}.

(a) Consider the manifold from \ref{case1} and the family
\begin{eqnarray*}
\omega_t=-{t\over 2}(t-3)\alpha_{14}+{1\over 4}(t^2-5t+4)(\alpha_{15}+
\alpha_{24})-{t\over 2}(t-3)(\alpha_{26}-
\alpha_{34})+\alpha_{16}-\alpha_{35}
\end{eqnarray*}
of closed 2-forms on $M$.
The form $\omega_t$ is symplectic if
$$
t^6-11t^5+39t^4-45t^3+4 t^2-20t+16\not=0.
$$

This polynomial has two real roots, and both of them lie out of the
interval $(1-\eps,4)$. So, because of \ref{example1}
$$
 h_4(\omega_2)=2,\quad h_4(\omega_1)=h_{4}(\omega_{3+\sqrt{17}\over 2})=3,
$$
and $h_4(\omega_t)=4$ for all other $t\in (1-\eps, 4)$.

(b) Consider the manifold from \ref{case2} and the family
$$
\omega_t=(1-t)\alpha_{13}-t(\alpha_{16}+
\alpha_{25}-\alpha_{34})+(1-t)(\alpha_{26}-\alpha_{45})
$$
of closed 2-forms. Since the polynomial $AE^2+BDE-CDE-D^3=3t^2-3t+1$ has no real
roots, we conclude that $\omega_t$ is a family of symplectic structures and
\begin{enumerate}
\item $h_4(\omega_1)=3$, $h_5(\omega_1)=0$;
\item $h_4(\omega_t)=4$, $h_5(\omega_t)=2$ for
$t\not= 1$.
\end{enumerate}
}
\end{comment}

\begin{remark}
{\rm Minding Theorem \ref{nil}, it is
natural to ask whether there exist flexible nilmanifolds of dimension
greater than 6. Taking into account the results of the next section,
we see that the answer is affirmative. However, the question remains
open for irreducible compact nilmanifolds of dimension greater than 6.
}
\end{remark}

\begin{remark}\label{fls}
{\rm We have also considered the 6-dimensional compact completely
solvable manifold $M$ constructed by Fern\'andez--de Le\'on--Saralegi
~\cite{FLS}. This compact manifold does not satisfy the Hard Lefschetz
condition (although it is of Lefschetz type). We have
$$
b_1(M)=b_5(M)=2, \ b_2(M)=b_4(M)=3, \ b_3(M)=4.
$$
Furthermore, $h_3(M,\omega)=4$ and $h_4(M,\omega)=h_5(M,\omega)=2$ for
every symplectic form $\omega$ on $M$. We do not explain the details
because $M$ is not flexible. }
\end{remark}

The following question also seems to be interesting.

{\bf Question}: Is there a Nomizu's type result for compact
nilmanifolds (more generally, homogeneous spaces) and the
symplectically harmonic cohomology $H^*_{\hr}(M)$. In another words, does
a symplectically harmonic de Rham cohomology class contain a
homogeneous symplectically harmonic representative (if we are
considering homogeneous symplectic structures)?

The answer is affirmative for degrees $k\leq 2$ and $k\geq 2m-2$.
Indeed, let $G/\Gamma$ be a $2m$-dimensional compact nilmanifold with
a homogeneous symplectic form $\omega$, and let $\g$ be the Lie algebra
of $G$. Since the image of a homogeneous form under each of the
operators $\ast$, $d$ and $L$ is homogeneous, the (finite dimensional)
subspaces $\Omega^{\ast}_{\hr}(\g^{\ast})$ and $\Lambda^{\ast}(\g^{\ast})$
are ${\frak s\frak l}(2)$-submodules of $\Omega^{\ast}(G/\Gamma)$. Therefore,
obvious analogs of Proposition \ref{dual}, Corollaries \ref{yan} and
\ref{tool} hold for $\g^{\ast}$ and the result follows from the
Nomizu's theorem.

\section{Product formula for symplectically harmonic cohomology}
\setcounter{equation}{0}

Let $(M^{2m},\omega_1)$ and $(N^{2n},\omega_2)$ be two symplectic
manifolds. Consider the symplectic product manifold $(M\times N, \omega)$ where
$\omega=p_1^{\#}\omega_1+p_2^{\#}\omega_2$ and
$$
p_1: M \times N \to M, \quad p_2: M \times N \to N
$$
are the projections. Given two forms $\alpha\in
\Omega^p(M)$ and $\beta\in \Omega^q(N)$, consider the form
$$
\alpha \boxtimes \beta :=(p_1^{\#}\alpha)\wedge (p_2^{\#}\beta)\in
\Omega^{p+q}(M\times N).
$$

\begin{proposition}\label{star}
{\rm(\cite{Br})}
$$
\ast(\alpha \boxtimes \beta)=(-1)^{pq}(\ast_{1}\alpha)\boxtimes
(\ast_{2}\beta),
$$
\hqed
\end{proposition}

\begin{corollary}\label{product}
{\rm (i)} $i(\Pi)(\alpha\boxtimes \beta)=(i(\Pi_{1})\alpha)\boxtimes
\beta+ \alpha\boxtimes (i(\Pi_{2})\beta)$;
\par {\rm (ii)} $\delta (\alpha\boxtimes \beta)=(\delta_1\alpha)\boxtimes
\beta+ (-1)^p \alpha\boxtimes (\delta_{2}\beta)$;
\par {\rm (iii)} $\Omega^p_{\hr}(M)\boxtimes \Omega^q_{\hr}(N)\subset \Omega^{p+q}_{\hr}(M\times N)$;
\par {\rm (iv)} for all $k$ we have
$$
\sum_{p+q=k} h_p(M)
h_{q}(N)\leq h_{k}(M\times N).
$$
\hqed
\end{corollary}

{\bf Question:} When the inequality in \ref{product}(iv) turns out to be the
equality?

\medskip Now we consider the Lefschetz map
$L_{[\omega]}^{m+n-k}:H^k(M\times N)\longrightarrow
H^{2m+2n-k}(M\times N)$.
\begin{proposition}\label{betti}
\begin{eqnarray*}
h_{2(m+n)-1}(M\times N)&=& h_{2m-1}(M) + h_{2n-1}(N),\\
h_{2(m+n)-2}(M\times N)&=&h_{2m-2}(M) + h_{2m-1}(M)h_{2n-1}(N) +
h_{2n-2}(N).
\end{eqnarray*}
\end{proposition}

\begin{proof}
Because of the K\"unneth isomorphism $H^*(M \times N) \cong
H^*(M)\otimes H^*(N)$, we conclude that
$$
\Im\, \left(L_{[\omega]}^{m+n-1}\right)=
[\omega_1]^m\otimes \left(\Im\, L_{[\omega_2]}^{n-1}\right)\oplus
\left(\Im\, L_{[\omega_1]}^{m-1}\right)\otimes [\omega_2]^n,
$$
and the first equality follows from \ref{tool}. Similarly,
\begin{eqnarray*}
\Im\, \left(L_{[\omega]}^{m+n-2}\right)&=&
\{([\omega_1]^{m-1}u)\otimes [\omega_2]^{n-1}+
([\omega_1]^{m-2} u)\otimes [\omega_2]^{n}\, |\, u\in H^2(M)\}\\
&\oplus &
\{[\omega_1]^{m}\otimes ([\omega_2]^{n-2}v)+
[\omega_1]^{m-1}\otimes ([\omega_2]^{n-1} v)\, |\, v\in H^2(N)\}\\
&\oplus & \{([\omega_1]^{m-1}w_1)\otimes ([\omega_2]^{n-1} w_2)\, |\,
w_1\in H^2(M),\,w_2 \in H^2(N)\}.
\end{eqnarray*}
Now, in view of \ref{tool}, the computation of dimensions completes
the proof.
\end{proof}

\begin{corollary}
Let $M^{2m}$ be a manifold which admits a family of symplectic forms
such that the symplectically harmonic Betti number $h_{2m-k}$ varies
for $k=1$ or $k=2$. Then $M \times N$ is a flexible manifold whenever
a manifold $N$ admits a symplectic structure.
\end{corollary}

\begin{proof}
If $h_{2m-1}(M)$ varies then the result follows from the first
equality of \ref{betti}. If $h_{2m-2}(M)$ varies but $h_{2m-1}(M)$
does not vary then the result follows from the second equality of
\ref{betti}.
\end{proof}

\section{Duality}
\setcounter{equation}{0}

Consider a symplectic manifold $(M^{2m},\omega)$ and the chain complex
$$
\cdots \hfl{}{}\, \Omega^{k+1}(M)\, \hfl{\delta}{} \,\Omega^{k}(M)
\,\hfl{\delta}{}\,\Omega^{k-1}(M)\,\hfl{}{} \cdots
$$
with $\delta$ as in (\ref{relation}). The following proposition
follows directly from the definition of $\delta$.

\begin{proposition}
\par{\rm (i)} $\delta \alpha =0$ if and only if $d *\alpha=0;$
\par{\rm (ii)} $\alpha \in
\Im\,\delta$ if and only if $*\alpha \in \Im\, d$.
\hqed
\end{proposition}

We define
$$
H^k_{\delta}(M)=H^k_{\delta}(M,\omega)=\Ker\,\delta^k/\Im\, \delta^{k+1}
\hbox{ where } \delta^i=\delta: \Omega^i \rightarrow \Omega^{i-1}.
$$

\begin{corollary}
The operator $*:\Omega^k \to \Omega^{2m-k}$ induces an isomorphism
$$
*:H^k(M) \to H^{2m-k}_{\delta}(M).
$$
In particular,
$H^k_{\delta}(M)=H^k(M)$ for $M$ closed.
\hqed
\end{corollary}

We dualize the definition of symplectically harmonic Betti numbers
$h_k$ by setting
$$
h^*_k(M)=h^*_k(M,\omega):=\dim\,\left(\Omega^k_{\hr}/\Im\,\delta
\cap\Omega^k_{\hr}\right)
$$

\begin{corollary}
$h^*_{m-k}(M)=h_{m+k}(M)$.
\hqed
\end{corollary}

In particular, in view of \ref{yan}, if $M$ is closed then
$h^*_{2m-k}(M)=b_{2m-k}(M)$ for $k=0,1,2$.

It is clear that many other results of Sections 2 and 3 can be
dualized in a similar way. We leave it to the reader.

\section*{Acknowledgments} This work was done during the stay of three of the
authors in the Mathematisches Forschungsinstitut in Oberwolfach under
the ``Research in Pairs'' program financed by the
``Volkswagen-Stiftung''. We express our sincere thanks to the
Institute and to the Foundation for the hospitality and wonderful
atmosphere. The third author was partially financed by the Polish
Research Commitee (KBN). The first and fourth authors were partially
supported by the project UPV 127.310-EA147/98 and DGICYT
PB97-0504-C02-02.

We thank Robert Gompf, Boris Khesin and Dusa McDuff for helpful discussions.

\end{document}